\documentclass[preprint,fleqn,twoside,10pt]{article}
\usepackage{amssymb}

\usepackage{amssymb,amsmath}
\usepackage{graphics}
\usepackage{graphicx}
\usepackage{color}
\newcommand{\dl}[1]{{\bf Theorem{#1.}}}
\newcommand{\yl}[1]{{\bf Lemma{#1.}}}

\newcommand{\de}[1]{{\bf Definition{#1}}}
\newcommand{\zm}{{\bf Proof.}}
\newcommand{\zb}{\textbf{$\hfill\Box$}}
\newcommand{\la}[1]{\label{#1}}
\newcommand{\rf}[1]{(\ref{#1})}

\newcommand{\commentout}[1]{{}}

\begin{document}
\thispagestyle{empty} \setcounter{page}{1}
{\small\noindent Fractals\\ \copyright World Scientific Publishing Company}\\

\begin{center}
{\bf A DIFFERENCE METHOD FOR SOLVING THE NONLINEAR $q$-FACTIONAL DIFFERENTIAL EQUATIONS ON TIME SCALES}\\[0.2cm]

TIE ZHANG\footnote{Corresponding author} and CAN TONG\\
{\small \em Department of Mathematics and the State Key Laboratory\\ of Synthetical Automation for Process Industries\\ Northeastern University, Shenyang 110004, China\\
$^*$ztmath@163.com}\\[1.5cm]

\begin{minipage}{5in}
\small\baselineskip 0.5cm
The $q$-fractional differential equation usually describe the physics process imposed on the time scale set $T_q$.
In this paper, we first propose a difference formula for discretizing the fractional $q$-derivative $^cD_q^\alpha x(t)$ on the time scale set $T_q$ with order $0<\alpha<1$ and scale index $0<q<1$. We establish a rigours truncation error boundness and prove that this difference formula is unconditionally stable. Then, we consider the difference method for solving the initial problem of $q$-fractional differential equation: $^cD_q^\alpha x(t)=f(t,x(t))$ on the time scale set. We prove the unique existence and stability of the difference solution and give the convergence analysis. Numerical experiments show the effectiveness  and high accuracy of the proposed difference method.\\

{\em Keywords:} The fractional $q$-derivative; Difference formula; Truncation error; Unconditional stability; The $q$-fractional differential equation; Convergence analysis.\\[0.5cm]
\end{minipage}
\end{center}

\section*{1. INTRODUCTION}
In recent years, the $q$-calculus (also called quantum calculus) and $q$-fractional calculus have attracted much attention in mathematical and physical fields, such as number theory, special functions, basic hyper-geometric functions, operator theory, quantum dynamics and stochastic analysis \cite{Annaby,Aral,Atici,Kac,Ezeafulukwe,Jarad,Umarov}.
The $q$-calculus is a class of calculus defined on a $q$-geometry set $T_q$ where $q$ is the scale index of set $T_q$. The $q$-fractional differential equations usually describe some physics processes imposed on one time scale set $T_{q,b}=\{q^nb, n=0,1,\cdots\}\cup\{0\}$ where $q$ is used to indicate the discrete time path on which the corresponding physical process varies. The $q$-derivative was originally presented by Jackson \cite{Jackson} in 1908 and the $q$-integral originated from a $q$-analogue of the usual integral proposed by Al-Salam \cite{Salam} in 1966. Following Al-Salam's work, Agarwal \cite{Agarwal} further gave the definition of $q$-fractional calculus. We refer the readers to \cite{Annaby,Aral,Kac,Rabia,Thabet,Andrews} for some basic introduction on $q$-calculus and $q$-fractional calculus. Compared with the classical fractional calculus, the study of $q$-fractional calculus is still immature.

The existence of solutions for the Caputo type $q$-fractional initial value and boundary value problems has been studied by many researches, see \cite{Ahmad,Almeida,Aktu,Chen,Zhang1}. By using the Lyapunov's direct method, Jarad et. al. \cite{Jarad} studied the stability of Caputo type $q$-fractional non-autonomous systems. Wang et. al. \cite{Wang} investigated the existence of extremal solutions of the Caputo type $q$-fractional differential equation. Zhang et. al. \cite{Tang, Zhang1} gave new existence theorems for the initial value problem of $q$-fractional differential equation  under different conditions of source term functions. On the discrete or numerical methods, Abdeljawad and Baleanu \cite{Abdeljawad} first presented the successive approximation method to obtain an approximate solution of the Caputo type $q$-fractional differential equation. Then, Salahshour et. al. \cite{Salahshour} investigated the convergence of the successive approximation method presented by Abdeljawad. Wu and Baleanu \cite{Wu} used the variational iteration method and the Lagrange multiplier method to study the Caputo type $q$-fractional initial value problem. Recently, Zhang and Tang \cite{Zhang} established a difference method for solving the initial problem of the $q$-fractional differential equation, but no stability and convergence analysis are given for this method. At present, only a few discrete or approximate methods are proposed for the $q$-fractional problems. In particular, almost no difference method is analyzed theoretically for solving the $q$-fractional differential equations.

In this article, we first present a difference formula (called the $L_{1,q}$ formula) to discretize the fractional $q$-derivative $^cD_q^\alpha x(t)$ of Caputo type on the time scale set $T_{q,b}=\{q^nb, n=0,1,\cdots\}\cup\{0\}$. Different from the standard difference formula, we must establish the difference formula on set $T_{q,b}$ , that is, the mesh points $\{t_k\}$ should be in set $T_{q,b}$. Thus, the mesh partition with nodes $\{t_k=bq^k,\,k=0,\cdots,N\}$ is nonuniform and the maximal mesh step does not vanish as $N\rightarrow \infty$. This makes the error analysis and stability estimate  much more difficult than that of the standard difference methods which can be established on any proper meshes, for example, the uniform mesh. Our difference formula is constructed by using the piecewise linear interpolation to approximate the integrand function on set $T_{q,b}$. The coefficients of this formula are analyzed in detail. We derive the truncation error boundness and prove that this difference formula is unconditionally stable. Then, using this difference formula, we establish a difference method for solving the initial problem of $q$-fractional differential equation: $^cD_q^\alpha x(t)=f(t,x(t)), t\in T_{q,b}$. We prove that this difference method is stale and give an error estimation of $\triangle t_n^2$-order. Numerical experiments show the high accuracy and effectiveness of this difference formula. To the authors' best knowledge, it is the first time that an unconditionally stable difference formula is presented and analyzed for the $q$-fractional problems. Our work provides a numerical approach for solving the $q$-fractional problems.

This paper is organized as follows. In Section 2, we introduce some basic notations and operations on $q$-fractional calculus. In Section 3, we present the $L_{1,q}$ difference formula and derive the truncation error boundness. Section 4 is contributed to the stability analysis of this difference formula. In Section 5, this difference formula is used to solve the initial value problem of $q$-fractional differential equation and the unique existence of the difference solution, stability and error estimation are given. Numerical experiments are provided in Section 6. Section 7 gives some conclusions.

\section*{2. THE $q$-CALCULUS AND $q$-FRACTIONAL CALCULUS}

In this section, we introduce some notations and operations on the $q$-calculus and the $q$-fractional calculus.

Let $\mathbb{N}=\left\{1, 2, 3, ...\right\}$ be the positive integer set. Define the $q$-shifted factorial operation by
\begin{equation}
(t-s)^{(0)}=1, (t-s)^{(k)}=\prod\limits_{i=0}^{k-1}(t-q^is),\,k\in \mathbb{N}.\la{2.1}
\end{equation}
If $\alpha\in R$ and $\pm\alpha\notin\mathbb{N}$, then
\begin{equation}
(t-s)^{(\alpha)}=t^{\alpha}\prod\limits_{i=0}^{\infty}\frac{t-q^is}{t-q^{\alpha+i}s},\,0\leq s\leq t.\la{2.2}
\end{equation}

Let $\mathbb{C}$ be the complex set. For $0<q<1$ and $\alpha\in\mathbb{C}\setminus \{-n: n\in \mathbb{N}\cup\{0\} \}$, the $q$-gamma function $\Gamma_q(\alpha)$ is defined by
\begin{equation}
\Gamma_q(\alpha)=(1-q)^{(\alpha-1)}(1-q)^{1-\alpha}.\la{2.3}
\end{equation}
Introduce the notations
$$
[\alpha]_q=\frac{1-q^\alpha}{1-q},\;\;[n]_q!=[n]_q[n-1]_q\cdots[1]_q.
$$
By a straightforward computation, it holds that
$$
\Gamma_q(1)=1,\;\Gamma_q(\alpha+1)=[\alpha]_q \Gamma_q(\alpha),\;\Gamma_q(n+1)=[n]_q!\,.
$$
For $\mathrm{Re}(\alpha)>0$, $\mathrm{Re}(\beta)>0$, the $q$-beta function is defined by
\begin{equation}
B_q(\alpha, \beta)=\int_{0}^{1}t^{\alpha-1}(1-qt)^{(\beta-1)}d_qt,\la{2.4}
\end{equation}
where the $q$-integral $\int \,d_qt$ is defined by \rf{2.5}--\rf{2.6}. The relation between the $q$-beta function and the $q$-gamma function is given as follows \cite{Askey}.
$$
B_q(\alpha, \beta)=\frac{\Gamma_q(\alpha)\Gamma_q(\beta)}{\Gamma_q(\alpha+\beta)}.
$$

For a fixed point $q\in R$, a subset $A\subset R$ is called $q$-geometric if $qt\in A$ whenever $t\in A$.
It is easy to see that $\forall\, t\in A$, $A$ contains all geometric sequences $\{t q^n\}_{n=0}^\infty$. A typical $q$-geometry set is the time scale set $T_q$ defined by $T_q=\left\{q^n: n\in \mathbb{Z}\right\}\cup \left\{0\right\}$, where $0<q<1$, $\mathbb{Z}=\left\{0, \pm1, \pm2,...\right\}$.\\
\de{ 1.} \cite{Jackson}\quad {\em Let $f(t)$ be a real valued function on set $A$, $0<q< 1$. The $q$-derivative of $f(t)$ is defined by}
\begin{eqnarray*}
D_q f(t)&=&\frac{d_qf(t)}{d_qt}=\frac{f(qt)-f(t)}{(q-1)t},\,\, t\in A\backslash \left\{0\right\},\\
D_q f(0)&=&\frac{d_qf(t)}{d_qt}|_{t=0}=\lim_{n\to\infty}\frac{f(tq^n)- f(0)}{t q^n},\,\, 0<q<1,\,t\neq 0.
\end{eqnarray*}

The $q$-derivative is a discrete analogy of the classical derivative. Obviously, for a function $f(t)$ defined on $R$, its $q$-derivative $D_qf(t)$ always exists for $t\not= 0$ and $D_qf(0)$ is well defined if the classical derivative $f'(0)$ exists.

The high order $q$-derivative $D_q^nf(t)$ is defined by $D_q^nf(t)=D_q(D^{n-1}_qf(t)),\,n\geq 2$.

Let $f(t)$ and $g(t)$ be real valued functions on $R$ and both $f(t)$ and $g(t)$ are differentiable at $t=0$. Then, the following $q$-derivative operation rules hold \cite{Annaby}
\begin{eqnarray*}
&&D_q\left(a f(t)\pm b g(t)\right)=a D_q f(t)\pm b D_q g(t),\,a,b\in R,\\
&&D_q\left(f(t) g(t)\right)=g(t)D_q f(t)+f(qt) D_q g(t),\\
&&D_q\big(\frac{f(t)}{g(t)}\big)=\frac{g(t)D_q f(t)-f(t) D_q g(t)}{g(t) g(qt)},\,g(t)\not=0,\,g(qt)\not=0.
\end{eqnarray*}
\de{ 2.} \cite{Jackson1}\quad {\em The $q$-integral on interval $(a,b)$ is defined by
\begin{equation}
\int_{a}^{b}f(t)d_q t=\int_{0}^{b}f(t)d_q t-\int_{0}^{a}f(t)d_q t,\,\, a,\,b\in A,\la{2.5}
\end{equation}
where}
\begin{equation}
\int_{0}^{x}f(t)d_q t=(1-q)\sum_{n=0}^{\infty}xq^n f(x q^n),\,\, x\in A.\la{2.6}
\end{equation}

By this definition, it is easy to see that
\begin{eqnarray}
&&\Big|\int_0^bf(t)d_qt\Big|\leq \int_0^b|f(t)|d_qt,\;0<b.\la{2.6a}\\
&&\int_a^bf(t)d_qt=\int_a^cf(t)d_qt+\int_c^bf(t)d_qt,\;a<c<b.\la{2.6b}
\end{eqnarray}

The following lemma gives the $q$-integration by parts rule.\\
\yl{ 1}\quad{\em Let $f(t)$ and $g(t)$ are real valued functions on set $A$, $0\leq a<b,\,a,b\in A$ and $0<q<1$, then}
\begin{eqnarray}
\int_{a}^{b}g(t) D_q f(t)d_q t=
(fg)(b)-(fg)(a)-\int_{a}^{b}f(q t) D_q g(t) d_q t,\la{2.7}\\
\int_{a}^{b}g(q t) D_q f(t)d_q t=(fg)(b)-(fg)(a)-\int_{a}^{b} f(t)D_q g(t) d_q t.\la{2.8}
\end{eqnarray}
\zm\quad We only prove conclusion \rf{2.7}, the proof of conclusion \rf{2.8} is completely similar. When $a=0$, we have from \rf{2.6} that
\begin{align*}
\int_{0}^{b}g(t)D_q f(t)d_q t&=(1-q)\sum_{n=0}^{\infty}b q^n g(b q^n) D_qf(b q^n)\\
&=(1-q)b\sum_{n=0}^{\infty}q^n g(bq^n)\frac{f(b q^n)-f(bq^{n+1})}{(1-q)b q^n}\\
&=\sum_{n=0}^{\infty}g(b q^n)\left(f(b q^n)-f(b q^{n+1})\right),
\end{align*}
and
\begin{align*}
\int_{0}^{b}f(qt)D_q g(t)d_q t&=(1-q)\sum_{n=0}^{\infty}b q^n f(b q^{n+1}) D_qg(b q^n)\\
&=(1-q)b\sum_{n=0}^{\infty}q^n f(bq^{n+1})\frac{g(b q^n)-g(bq^{n+1})}{(1-q)b q^n}\\
&=\sum_{n=0}^{\infty}f(b q^{n+1})\left(g(b q^n)-g(b q^{n+1})\right).
\end{align*}
So
\begin{eqnarray}
&&\int_{0}^{b}g(t)D_q f(t)d_q t+\int_{0}^{b}f(qt)D_q g(t)d_q t
=\sum_{n=0}^{\infty}\Big(f(bq^n)g(bq^n)-f(bq^{n+1})g(bq^{n+1})\Big)\nonumber\\
&&=f(b)g(b)-\lim_{n\to\infty}(fg)(bq^{n+1})=(fg)(b)-(fg)(0).\la{2.9}
\end{eqnarray}
 When $a>0$, we have from \rf{2.5} and \rf{2.9}
\begin{eqnarray*}
\int_{a}^{b}g(t)D_q f(t)d_q t&=&\int_{0}^{b}g(t)D_q f(t)d_q t-\int_{0}^{a}g(t)D_q f(t)d_q t\\
&=&(fg)(b)-(fg)(0)-\int_{0}^{b}f(qt) D_q g(t)d_q t\\
&&-(fg)(a)+(fg)(0)+\int_{0}^{a}f(qt) D_q g(t)d_q t\\
&=&(fg)(b)-(fg)(a)-\int_{a}^{b}f(qt) D_q g(t)d_q t.
\end{eqnarray*}
The proof is completed. \zb

Now we introduce the $q$-fractional calculus.\\
\de{ 3.} \cite{Atici}\quad{\em Let $t\in A$, $\alpha \neq -1, -2, ...$ and $a\geq 0$. The $\alpha$-order $q$-fractional integral of the Riemann-Liouville type with the lower limit point $a$ is defined by $I_{q, a}^0 f(t)=f(t)$ and}
\begin{equation}
I_{q,a}^\alpha f(t)=\frac{1}{\Gamma_q(\alpha)}\int_{a}^{t}(t-qs)^{(\alpha-1)}f(s)d_qs.\la{2.10}
\end{equation}
\de{ 4.} \cite{Rajkovi}\quad{\em Let $a\geq 0, a\in A, n=\lceil\alpha\rceil$. The $\alpha$-order Riemann-Liouville type fractional $q$-derivative of  function $f(t): (a,\infty) \rightarrow \mathrm{R}$ is defined by
\begin{eqnarray}
D_{q,a}^\alpha f(t)=\left\{
\begin{array}{lll}
&I_{q,a}^{-\alpha}f(t),\;&\alpha \leq 0,\\
&D_{q}^nI_{q,a}^{n-\alpha}f(t),&\alpha >0,
\end{array}\right.\la{2.11}
\end{eqnarray}
where $\lceil\alpha\rceil$ denotes the smallest integer that is greater or equal to $\alpha$.}\\
\de{ 5.} \cite{Rajkovi}\quad{\em Let $a\geq 0,\,a\in A,\,n=\lceil\alpha\rceil$. The $\alpha$-order Caputo type fractional $q$-derivative of function $f(t): (a,\infty) \rightarrow \mathrm{R}$ is defined by}
\begin{eqnarray}
^{c}D_{q,a}^\alpha f(t)=\left\{
\begin{array}{lll}
&I_{q,a}^{-\alpha}f(t), \;&\alpha \leq 0,\\
&I_{q,a}^{n-\alpha} D_{q}^nf(t),&\alpha >0.
\end{array}\right.\la{2.12}
\end{eqnarray}

For simplicity, we often use the notations $I_q^\alpha f(t)$ instead of $I_{q, 0}^\alpha f(t)$, $D_{q}^\alpha f(t)$ instead of $D_{q,0}^\alpha f(t)$ and $^{c}D_{q}^\alpha f(t)$ instead of $^{c}D_{q,0}^\alpha f(t)$, respectively.

Under a certain condition, the following relations between the Caputo type fractional $q$-derivative and the Riemann-Liouville type fractional $q$-derivative hold \cite{Annaby,Abdeljawad}
\begin{equation}
^cD_q^\alpha f(t)=D_q^\alpha\Big(f(t)-\sum_{j=0}^{n-1}\frac{D_q^j f(0)}{\Gamma_q(j+1)}t^j\Big),\;n=\lceil \alpha\rceil,\,\alpha>0.\la{2.13}
\end{equation}
In particular, for $0<\alpha<1$,
\begin{equation}
^cD_{q}^\alpha f(t)=D_{q}^\alpha f(t)-\frac{t^{-\alpha}}{\Gamma_q(1-\alpha)}f(0).\la{2.14}
\end{equation}

\section*{3. THE $L_{1,q}$ DIFFERENCE FORMULA}
In this section, we establish a difference formula (called the $L_{1,q}$ formula) for discretizing the Caputo type fractional $q$-derivative $^cD_{q}^\alpha x(t)$ with $0<\alpha,q<1$ and give the truncation error boundness.

Let us consider the Caputo type fractional $q$-derivative:
\begin{equation}
^cD_{q}^\alpha x(t)=\frac{1}{\Gamma_q(1-\alpha)}\int_0^t(t-qs)^{(-\alpha)}D_qx(s)d_qs,\,0<\alpha,q<1,\;0<t\leq b,\,t\in T_{q,b}.\la{3.1}
\end{equation}
We need to construct a difference formula to approximate $^cD_{q}^\alpha x(t)$ on the discrete point set $\{t_k\}\subset T_{q,b}$ where the time scale set $T_{q,b}=\{bq^n: n=0,1,\dots\}\bigcup\{0\}$.

Let $0=t_0< t_1<\cdots < t_N=b$ be a nonuniform partition of $[0,b]$ with mesh point $t_k=bq^{N-k}\in T_{q,b}$ and step sizes $\Delta t_k=t_k-t_{k-1},\,1\leq k\leq N$ where $N\ge 1$ is a positive integer. Introduce the piecewise linear interpolation of function $x(s)$:
\begin{equation}
L_{1,k}(s)=\frac{s-t_{k-1}}{\Delta t_k} x(t_k)+\frac{t_{k}-s}{\Delta t_k} x(t_{k-1}), \,s\in [t_{k-1}, t_k],\, k=1, 2, ..., N,\la{3.2}
\end{equation}
with the interpolation error
\begin{equation}
R_k(s)=x(s)-L_{1,k}(s), s\in [t_{k-1},t_k],\;k=1, 2, ..., N.\la{3.3}
\end{equation}
Since
\begin{equation}
^cD_{q}^\alpha x(t_n)=\frac{1}{\Gamma_q(1-\alpha)}\sum_{k=1}^n\int_{t_{k-1}}^{t_k}(t_n-qs)^{(-\alpha)}D_qx(s)d_qs,\la{3.4}
\end{equation}
replacing $D_qx(s)$ by $D_qL_{1,k}x(s)$ in \rf{3.4} and noting that $D_qL_{1,k}x(s)=(x(t_k)-x(t_{k-1}))/\triangle t_k$, we obtain
\begin{eqnarray}
&&^cD_{q}^\alpha x(t_n)=\frac{1}{\Gamma_q(1-\alpha)}\sum_{k=1}^nb_k\big(x(t_k)-x(t_{k-1})\big)+R^n_q,\la{3.5}\\
&&b_k=\frac{1}{\triangle t_k}\int_{t_{k-1}}^{t_k}(t_n-qs)^{(-\alpha)}d_qs,\;\,k=1,2,\dots, n\,,\la{3.6}
\end{eqnarray}
where the truncation error
\begin{equation}
R^n_q=\frac{1}{\Gamma_q(1-\alpha)}\sum_{k=1}^n\int_{t_{k-1}}^{t_k}(t_n-qs)^{(-\alpha)}D_qR_k(s)d_qs.\la{3.7}
\end{equation}
Thus, we derive the $L_{1,q}$ difference formula:
\begin{equation}
\triangle_q^\alpha x^n
\doteq\frac{1}{\Gamma_q(1-\alpha)}\sum_{k=1}^nb_k\big(x^k-x^{k-1}),\;\;^cD_{q}^\alpha x(t_n)=\triangle ^\alpha_qx(t_n)+R^n_q\,,\la{3.8}
\end{equation}
where $x^k$ is any one mesh function defined on $\{t_k\}$ and when $x(t)$ is continuous on $[0,t_N]$, we set $x^k=x(t_k)$.\\
{\bf Remark 1.}\quad In constructing the difference scheme, we used $D_qL_{1,k}(t)$ to approximate the integrand function $D_qx(t)$ on interval $[t_{k-1},t_k]$ in (20). Since (noting that the $q$-mesh point $t_{k-1}=qt_k$)
$$
D_qL_{1,k}(t)=\frac{x(t_k)-x(t_{k-1})}{t_k-t_{k-1}}=\frac{x(t_k)-x(qt_{k})}{(1-q)t_k}=D_qx(t_k),\;t\in [t_{k-1},t_k],
$$
so actually, we used the $q$-rectangle quadrature formula to discretize the integral on $[t_{k-1},t_k]$.

Below we estimate the truncation error $R^n_q$. We first give several lemmas.\\
\yl{ 2}($q$-Rolle lemma \cite{Rajkovi1})\quad{\em Let $x(t)$ be a continuous function on $[a,b]$ and $x(a)=x(b)$. Then, for $q \in (0,1)$, there exists point $\xi\in(a,b)$ such that $D_q x(\xi)=0$.}

By means of Lemma 2, we can establish the following interpolation error formula.\\
\yl{ 3}\quad{\em Assume that $x(s)$ and $D_q x(s)$ is continuous on $[t_{k-1},t_k]$. Then, the linear interpolation error function $R_k(s)=x(s)-L_{1,k}(s)$ has the following expression}
$$
R_k(s)=\frac{1}{1+q}D_q^2 x(\xi_k)(s-t_{k-1})(s-t_k), s\in[t_{k-1},t_k], \,\xi_k \in (t_{k-1},t_k),\, 1\le k\le N.
$$
\zm\quad Since $R_k(t_{k-1})=R_k(t_k)=0$, we may assume that
\begin{equation}
R_k(s)=K_q(s)(s-t_{k-1})(s-t_{k}),\, s\in [t_{k-1},t_{k}],\la{3.9}
\end{equation}
where $K_q(s)$ is an undetermined function. For any fixed point $s\in (t_{k-1},t_{k})$, let $\varphi_q(t)=x(t)-L_{1,k}(t)-K_q(s)(t-t_{k-1})(t-t_k)$. Since $\varphi_q(t_{k-1})=\varphi_q(t_k)=\varphi_q(s)=0$, by using the $q$-Rolle lemma twice, we know that $D_q^2\varphi_q(t)$ has at least one zero point $\xi_k\in (t_{k-1},t_k)$.  So
\begin{eqnarray*}
D_q^2\varphi_q(\xi_k)&=&D^2_q\Big(x(t)-L_{1,k}(t)-K_q(s)(t-t_{k-1})(t-t_k)\Big)\Big|_{t=\xi_k}\\
&=&D_q^2x(\xi_k)-(1+q) K_q(s)=0.
\end{eqnarray*}
Hence, we obtain $ K_q(s)=D_q^2x(\xi_k)/(1+q)$. Substituting $K_q(s)$ into \rf{3.9}, the proof is completed.\zb\\
\yl{ 4}\quad{\em Let $0<\alpha, q<1$ and $D_q$ be the $q$-derivative operator with respect to variable $s$. Then}
\begin{eqnarray}
&&D_q(t-s)^{(-\alpha)}=-[-\alpha]_q(t-qs)^{(-\alpha-1)},\,0\leq s\leq t,\la{3.10}\\
&&|(t-qs)^{(-\alpha-1)}|\leq t^{-\alpha-1}\frac{1}{1-q^{\alpha}}\frac{1}{1-q^{1-\alpha}},\;0\leq s\leq t.\la{3.11}
\end{eqnarray}
\zm\quad We first prove conclusion \rf{3.10}. By the definition of operator $D_q$ and \rf{2.2}, we have
\begin{eqnarray}
&&D_q(t-s)^{(-\alpha)}=\frac{(t-qs)^{(-\alpha)}-(t-s)^{(-\alpha)}}{(q-1)s}\nonumber\\
&=&\frac{t^{-\alpha}}{(q-1)s}\lim_{m\rightarrow \infty}S_m,
\;S_m=\prod_{i=0}^m\frac{(t-q^{i+1}s)}{(t-q^{i+1-\alpha}s)}-\prod_{i=0}^m\frac{(t-q^{i}s)}{(t-q^{i-\alpha}s)}.\la{3.12}
\end{eqnarray}
By a straightforward computation, it yields
\begin{eqnarray*}
S_m&=&\prod_{i=1}^m\frac{(t-q^{i}s)}{(t-q^{i-\alpha}s)}\Big[\frac{(t-q^{m+1}s)}{(t-q^{m+1-\alpha}s)}
-\frac{(t-s)}{(t-q^{-\alpha})s}\Big]\\
&=&\prod_{i=1}^m\frac{(t-q^{i}s)}{(t-q^{i-\alpha}s)}\Big[\frac{st(1-q^{-\alpha})(1-q^{m+1})}{(t-q^{m+1-\alpha}s)(t-q^{-\alpha}s)}\Big]\\
&=&\prod_{i=0}^m\frac{(t-q^{i+1}s)}{(t-q^{i-\alpha}s)}\Big[\frac{st(1-q^{-\alpha})(1-q^{m+1})}{(t-q^{m+1-\alpha}s)(t-q^{m+1}s)}\Big]\\
&=&\prod_{i=0}^\infty\frac{(t-q^{i+1}s)}{(t-q^{i-\alpha}s)}\Big[\frac{s(1-q^{-\alpha})}{t}\Big],\;m\rightarrow \infty.
\end{eqnarray*}
Substituting this into \rf{3.12}, we obtain
$$
D_q(t-s)^{(-\alpha)}=\frac{(1-q^{-\alpha})}{(q-1)}t^{-\alpha-1}\prod_{i=0}^\infty\frac{(t-q^{i+1}s)}{(t-q^{i-\alpha}s)}
=-[-\alpha]_q(t-qs)^{(-\alpha-1)}.
$$

Next, we consider estimation \rf{3.11}. Since
\begin{eqnarray}
(t-qs)^{(-\alpha-1)}=t^{-\alpha-1}\lim_{m\rightarrow \infty}S_m',\;S'_m=\prod_{i=0}^m\frac{(t-q^{i+1}s)}{(t-q^{i-\alpha}s)},\la{3.13}
\end{eqnarray}
and
\begin{eqnarray*}
&&\max_{0\leq s\leq t}\big|\frac{(t-qs)}{(t-q^{-\alpha}s)}\big|=\max\big\{1,\big|\frac{1-q}{1-q^{-\alpha}}\big|\big\}\leq \frac{1-q}{1-q^\alpha},\\
&&\max_{0\leq s\leq t}\frac{(t-q^{i+1}s)}{(t-q^{i-\alpha}s)}=\frac{(1-q^{i+1})}{(1-q^{i-\alpha})},\;i\geq 1,
\end{eqnarray*}
so
\begin{eqnarray*}
|S'_m|&\leq& \frac{1-q}{1-q^\alpha}\prod_{i=1}^m\frac{(1-q^{i+1})}{(1-q^{i-\alpha})}\\
&=&\frac{1}{1-q^\alpha}\frac{1}{1-q^{1-\alpha}}
\frac{1-q}{1-q^{2-\alpha}}\frac{1-q^2}{1-q^{3-\alpha}}\cdots\frac{1-q^{m-1}}{1-q^{m-\alpha}}\frac{1-q^{m}}{1}\frac{1-q^{m+1}}{1}\\
&\leq& \frac{1}{1-q^{\alpha}}\frac{1}{1-q^{1-\alpha}}.
\end{eqnarray*}
Substituting this into \rf{3.13}, the proof is completed.\zb

Now we can give the truncation error estimation.\\
\dl{ 1}\quad {\em Let $x(t)$ and $D_qx(t)$ are continuous on $[0,b]$. Then, the truncation error of the $L_{1,q}$ difference formula \rf{3.8} satisfies the following estimate}
\begin{equation}
|R_q^n|\leq \frac{1}{4\Gamma_q(1-\alpha)}\frac{1}{1-q^2}\frac{1}{q^\alpha-q}\,t_n^{-\alpha}\triangle t_n^2\max_{0\leq t\leq t_n}|D_q^2x(t)|,\;\,n\geq 1.\la{3.14}
\end{equation}
\zm\quad Let $\bar{R}(s)=R_k(s),\,s\in [t_{k-1}, t_k],\,1\leq k\leq N$. From \rf{3.7}, the $q$-integration by parts formula \rf{2.8} and Lemma 4, we obtain
\begin{eqnarray*}
R^n_q&=&
\frac{1}{\Gamma_q(1-\alpha)}\sum_{k=1}^n\int_{t_{k-1}}^{t_k}(t_n-qs)^{(-\alpha)}D_qR_k(s)d_qs\\
&=&\frac{(t_n-qs)^{(-\alpha)}}{\Gamma_q(1-\alpha)}R_k(s)\big|_{t_{k-1}}^{t^k}-\frac{1}{\Gamma_q(1-\alpha)}\sum_{k=1}^n\int_{t_{k-1}}^{t_k}D_q(t_n-s)^{(-\alpha)}R_k(s)d_qs\\
&=&\frac{[-\alpha]_q}{\Gamma_q(1-\alpha)}\int_{0}^{t_n}(t_n-qs)^{(-\alpha-1)}\bar{R}(s)d_qs.
\end{eqnarray*}
Hence, using \rf{2.6a}--\rf{2.6b}, Lemma 3 and \rf{3.11} it yields
\begin{eqnarray*}
|R_q^n|&\leq&\frac{|[-\alpha]_q|}{\Gamma_q(1-\alpha)}\int_{0}^{t_n}|(t_n-qs)^{(-\alpha-1)}\bar{R}(s)|d_qs\\
&=&\frac{|[-\alpha]_q|}{\Gamma_q(1-\alpha)}\sum_{k=1}^{n}\int_{t_{k-1}}^{t_k}|(t_n-qs)^{(-\alpha-1)}R_k(s)|d_qs\\
&\leq& \frac{|[-\alpha]_q|}{\Gamma_q(1-\alpha)}\frac{1}{1+q}\frac{1}{4}\max_{1\leq k\leq n}|\triangle t_k|^2\max_{0\leq t\leq t_n}|D_q^2x(t)|\int_0^{t_n}|(t_n-qs)^{(-\alpha-1)}|d_qs\\
&\leq& \frac{|[-\alpha]_q|}{\Gamma_q(1-\alpha)}\frac{1}{1+q}\frac{1}{4}\max_{1\leq k\leq n}|\triangle t_k|^2\max_{0\leq t\leq t_n}|D_q^2x(t)|t_n^{-\alpha}\frac{1}{1-q^{\alpha}}\frac{1}{1-q^{1-\alpha}}\\
&=&\frac{1}{\Gamma_q(1-\alpha)}\frac{q^{-\alpha}}{1-q^2}\frac{1}{4}\frac{1}{1-q^{1-\alpha}}t_n^{-\alpha}\triangle t_n^2\max_{0\leq t\leq t_n}|D_q^2x(t)|\,.
\end{eqnarray*}
This gives the desired estimate.\zb

\section*{4. STABILITY OF THE $L_{1,q}$ DIFFERENCE FORMULA}

We first give two useful lemmas.\\
\yl { 5}\quad{\em For $0<q,\alpha<1$, the following inequality holds}
\begin{equation}
t_n^{-\alpha}<(t_n-q^{i+1}s)^{(-\alpha)}\leq(t_n-qs)^{(-\alpha)},\;i\geq 0,\;0\leq s\leq t_n.\la{4.1}
\end{equation}
\zm\quad For the left side inequality, since $(t_n-q^{i+j+1}s)/(t_n-q^{i+j+1-\alpha}s)>1$ for any $i,j\geq 0$, we have
$$
(t_n-q^{i+1}s)^{(-\alpha)}=t_n^{-\alpha}\prod_{j=0}^\infty\frac{(t_n-q^{i+j+1}s)}{(t_n-q^{i+j+1-\alpha}s)}> t_n^{-\alpha},\,i\geq 0,\,0\leq s\leq t_n.
$$
Next, by a straightforward computation, it yields
\begin{eqnarray*}
&&(t_n-q^{i+1}s)^{(-\alpha)}-(t_n-qs)^{(-\alpha)}\\
&=&t_n^{-\alpha}\prod_{j=0}^\infty\frac{(t_n-q^{i+j+1}s)}{(t_n-q^{i+j+1-\alpha}s)}
-t_n^{-\alpha}\prod_{j=0}^\infty\frac{(t_n-q^{j+1}s)}{(t_n-q^{j+1-\alpha}s)}\\
&=&t_n^{-\alpha}\prod_{j=i}^\infty\frac{(t_n-q^{j+1}s)}{(t_n-q^{j+1-\alpha}s)}
\Big(1-\prod_{j=0}^{i-1}\frac{(t_n-q^{j+1}s)}{(t_n-q^{j+1-\alpha}s)}\Big)<0,\,i\geq 1.
\end{eqnarray*}
The proof is completed.\zb\\
\yl{ 6}\quad{\em Let coefficient series $\{b_k\}$ be defined by \rf{3.6}. Then, it holds}
\begin{eqnarray}
&&b_k=(t_n-qt_k)^{(-\alpha)},\;k\geq 2,\la{4.2}\\
&&t_n^{-\alpha}<b_1<b_2<\cdots<b_{k-1}<b_k,\;2\leq k\leq n\,.\la{4.3}
\end{eqnarray}
\zm\quad Since $\triangle t_k=t_k-t_{k-1}=t_k(1-q),\,t_{k-1}=qt_k,\,k\geq 2$, then by \rf{3.6} we have
\begin{eqnarray*}
b_k&=&\frac{1}{\triangle t_k}\int_{t_{k-1}}^{t_k}(t_n-qs)^{(-\alpha)}d_qs\\
&=&\frac{1}{\triangle t_k}\int_{0}^{t_k}(t_n-qs)^{(-\alpha)}d_qs-\frac{1}{\triangle t_k}\int^{t_{k-1}}_{0}(t_n-qs)^{(-\alpha)}d_qs\\
&=&(1-q)\sum_{i=0}^\infty\frac{t_k}{\triangle t_k}q^i(t_n-q^{i+1}t_k)^{(-\alpha)}-(1-q)\sum_{i=0}^\infty\frac{t_{k-1}}{\triangle t_k}q^i(t_n-q^{i+1}t_{k-1})^{(-\alpha)}\\
&=&\sum_{i=0}^\infty q^i(t_n-q^{i+1}t_k)^{(-\alpha)}-\sum_{i=0}^\infty q^{i+1}(t_n-q^{i+2}t_{k})^{(-\alpha)}=(t_n-qt_k)^{(-\alpha)}.
\end{eqnarray*}
This gives equality \rf{4.2}. For \rf{4.3}, we first obtain from Lemma 5
\begin{eqnarray}
b_1&=&\frac{1}{\triangle t_1}\int_0^{t_1}(t_n-qs)^{(-\alpha)}d_qs=(1-q)\sum_{i=0}^\infty\frac{t_1}{\triangle t_1}q^i(t_n-q^{i+1}t_1)^{(-\alpha)}\nonumber\\
&>&t_n^{-\alpha}(1-q)\sum_{i=0}^\infty q^i=t_n^{-\alpha},\la{4.4}
\end{eqnarray}
and
\begin{eqnarray}
b_1&=&(1-q)\sum_{i=0}^\infty\frac{t_1}{\triangle t_1}q^i(t_n-q^{i+1}t_1)^{(-\alpha)}\nonumber\\
&\leq&(1-q)\sum_{i=0}^\infty q^i(t_n-qt_1)^{(-\alpha)}=(t_n-qt_1)^{(-\alpha)}.\la{4.5}
\end{eqnarray}
Next, it follows from \rf{4.2}, \rf{4.5} and Lemma 5
\begin{equation}
b_2-b_1\geq(t_n-qt_2)^{(-\alpha)}-(t_n-qt_1)^{(-\alpha)}=(t_n-qt_2)^{(-\alpha)}-(t_n-q^2t_2)^{(-\alpha)}>0.\la{4.6}
\end{equation}
Moreover, from \rf{4.2} and Lemma 5, it yields
\begin{eqnarray}
b_k-b_{k-1}&=&(t_n-qt_k)^{(-\alpha)}-(t_n-qt_{k-1})^{(-\alpha)}\nonumber\\
&=&(t_n-qt_k)^{(-\alpha)}-(t_n-q^2t_{k})^{(-\alpha)}>0,\;k\geq 3.\la{4.7}
\end{eqnarray}
Combining \rf{4.4}, \rf{4.6} and \rf{4.7}, estimation \rf{4.3} is derived.\zb

Now let us consider the stability of the $L_{1,q}$ difference formula \rf{3.8}. For this end, we consider the following difference equation:
\begin{equation}
\triangle_q^\alpha x^n=f^n,\;n=1,2,\dots, N,\;\;\;\triangle_q^\alpha x^n\doteq\frac{1}{\Gamma_q(1-\alpha)}\sum_{k=1}^nb_k(x^k-x^{k-1}),\la{4.8}
\end{equation}
where the initial value $x^0$ and the source term $f^n$ are given. Using the identity
$$
\sum_{k=1}^nb_k(x^k-x^{k-1})=b_nx^n+\sum_{k=1}^{n-1}(b_k-b_{k+1})x^k-b_1x^0,
$$
difference equation \rf{4.8} can be rewritten as
\begin{equation}
b_nx^n=b_1x^0+\sum_{k=1}^{n-1}(b_{k+1}-b_{k})x^k+\Gamma_q(1-\alpha)f^n,\;n\geq 1\,.\la{4.9}
\end{equation}
\dl{ 2}\quad{\em The $L_{1,q}$ difference formula is unconditionally stable such that the solution of difference equation \rf{4.9} satisfies}
$$
|x^n|\leq |x^0|+\Gamma_q(1-\alpha)t_n^{\alpha}\max_{1\leq k\leq n}|f^k|,\;n\geq 1.
$$
\zm\quad It follows from \rf{4.9} and Lemma 6 that
\begin{eqnarray}
b_n|x^n|&\leq& b_1|x^0|+\sum_{k=1}^{n-1}(b_{k+1}-b_k)\max_{1\leq k\leq n}|x^k|+\Gamma_q(1-\alpha)|f^n|\nonumber\\
&\leq& b_1|x^0|+(b_n-b_1)\max_{1\leq k\leq n}|x^k|+\Gamma_q(1-\alpha)|f^n|\,.\la{4.10}
\end{eqnarray}
Let $1\leq n_0\leq n$ be such that $|x^{n_0}|=\displaystyle{\max_{1\leq k\leq n}}|x^k|$. Taking $n=n_0$ in \rf{4.10}, it yields
$$
b_1\max_{1\leq k\leq n}|x^n|\leq b_1|x^0|+\Gamma_q(1-\alpha)\max_{1\leq k\leq n}|f^k|.
$$
Hence, using \rf{4.3}, we complete the proof.\zb

\section*{5. DIFFERENCE METHOD FOR THE $q$-FRACTIONAL\\ DIFFERENTIAL EQUATION}
Consider the initial value problem of $q$-fractional differential equation:
\begin{eqnarray}
\left\{
\begin{array}{lll}
&^cD_q^\alpha x(t)=f(t,x(t)),\;0<t\leq b,\;0<q,\,\alpha<1,\;t\in T_{q,b},\\
&x(0)=x^0,
\end{array}\right.\la{5.1}
\end{eqnarray}
where $f:[0,b]\times R^d\rightarrow R^d$ is continuous, $d\geq 1$. Let the partition of interval $[0,b]$ be given as that in Section 3. By means of the $L_{1,q}$ difference formula, we establish the difference method for solving problem \rf{5.1} as follows (see \rf{4.9})
\begin{eqnarray}
b_nx^n=b_1x^0+\sum_{k=1}^{n-1}(b_{k+1}-b_k)x^k+\Gamma_q(1-\alpha)f(t_n,x^n),\;n=1,2,\dots, N.\la{5.2}
\end{eqnarray}

We first consider the unique existence of the difference solution. \\
\dl{ 3}\quad{\em Assume that $f:[0,b]\times R^d\rightarrow R^d$ is continuous and satisfies the Lipschitz condition:
\begin{equation}
|f(t,x)-f(t,x')|\leq L|x-x'|,\;x,x'\in R^d,\;L_1=L\Gamma_q(1-\alpha)b^\alpha< 1.\la{5.3}
\end{equation}
Then, difference equation \rf{5.2} has a unique solution.}\\
\zm\quad We only need to prove that for any fixed $n\geq 1$, when $\{x^k\}_{k=0}^{n-1}$ are known, there exists a unique $x^n$ satisfying Eq. \rf{5.2}. Consider the iteration scheme:
\begin{equation}
b_nx^{n,l}=b_1x^0+\sum_{k=1}^{n-1}(b_{k+1}-b_k)x^k+\Gamma_q(1-\alpha)f(t_n,x^{n,l-1}),\,l=1,2,\dots,\;x^{n,0}=x^{n-1}.\la{5.4}
\end{equation}
Noting that $b_n\geq t_n^{-\alpha}\geq b^{-\alpha}$, we obtain from \rf{5.3}--\rf{5.4} that
\begin{eqnarray*}
|x^{n,l+1}-x^{n,l}|&=& b_n^{-1}\Gamma_q(1-\alpha)|f(t_n,x^{n,l})-f(t_n,x^{n,l-1})|\\
&\leq& L_1|x^{n,l}-x^{n,l-1}|\leq \cdots\leq L_1^l|x^{n,1}-x^{n,0}|,\;l\geq 1.
\end{eqnarray*}
Since $L_1<1$, this implies that the series $\{x^{n,l}\}$ is convergent. Letting $l\rightarrow \infty$ in Eq. \rf{5.4}, we see that $x^n\doteq x^{n,\infty}$ is the solution of difference equation \rf{5.2}. Now for fixed $n$, assume that Eq.  \rf{5.2} has two solutions $x^n$ and $\bar{x}^n$. Then, $x^n$ and $\bar{x}^n$ satisfy
$$
|x^n-\bar{x}^n|= b_n^{-1}\Gamma_q(1-\alpha)|f(t_n,x^{n})-f(t_n,\bar{x}^n)|\leq L_1|x^n-\bar{x}^n|<|x^n-\bar{x}^n|\,.
$$
This implies $x^n=\bar{x}^n$, the uniqueness is proved.  \zb

Next, we consider the stability of the difference solution.\\
\dl{ 4}\quad{\em Assume that condition \rf{5.3} holds. Then, the solution of difference equation \rf{5.2} satisfies the stability estimation}
$$
|x^n|\leq \frac{1}{1-L_1}\big[\,|x^0|+\Gamma_q(1-\alpha)t_n^{\alpha}\max_{1\leq k\leq n}|f(t_k,0)|\,\big],\;n\geq 1\,.
$$
\zm\quad Let $x^n$ be the solution of Eq. \rf{5.2}. By a similar argument to that of Theorem 2, we have
\begin{equation}
|x^n|\leq |x^0|+\Gamma_q(1-\alpha)t_n^{\alpha}\max_{1\leq k\leq n}|f(t_k,x^k)|,\;n\geq 1.\la{5.5}
\end{equation}
Using condition \rf{5.3}, it yields
$$
|f(t_k,x^k)|=|f(t_k,x^k)-f(t_k,0)|+|f(t_k,0)|\leq L|x_k|+|f(t_k,0)|\,.
$$
Combining this with \rf{5.5}, we arrive at
\begin{equation}
|x^n|\leq |x^0|+L_1\max_{1\leq k\leq n}|x^k|+\Gamma_q(1-\alpha)t_n^{\alpha}\max_{1\leq k\leq n}|f(t_k,0)|,\,n\geq 1\,.\la{5.6}
\end{equation}
Now, let $1\leq n_0\leq n$ be such that $|x^{n_0}|=\displaystyle{\max_{1\leq k\leq n}}|x^k|$. Then, taking $n=n_0$ in \rf{5.6}, we can derive
$$
(1-L_1)\max_{1\leq k\leq n}|x^k|\leq |x^0|+\Gamma_q(1-\alpha)t_n^{\alpha}\max_{1\leq k\leq n}|f(t_k,0)|,\;n\geq 1\,.
$$
This gives the desired estimation.\zb

Finally, we can give the error estimation.\\
\dl{ 5}\quad{\em Let $x(t)$ and $x^n$ be the solution of Eqs. \rf{5.1} and \rf{5.2}, respectively. Assume that $x(t)$ and $D_qx(t)$ are continuous on $[0,b]$ and condition \rf{5.3} holds. Then, we have the following error estimation.}
\begin{equation}
|x(t_n)-x^n|\leq \frac{1}{1-L_1}\frac{1}{4}\frac{1}{1-q^2}\frac{1}{q^\alpha-q}\triangle t_n^2\max_{0\leq t\leq t_n}|D_q^2x(t)|,\;n\geq 1.\la{5.7}
\end{equation}
\zm\quad From Eq. \rf{5.1}, we see that the solution $x(t)$ satisfies the discrete equation: $\triangle_q^\alpha x(t_n)=f(t_n,x(t_n))-R_q^n$, or
\begin{equation}
b_nx(t_n)=b_1x^0+\sum_{k=1}^{n-1}(b_{k+1}-b_k)x(t_k)+\Gamma_q(1-\alpha)f(t_n,x(t_n))-\Gamma_q(1-\alpha)R_q^n,\la{5.8}
\end{equation}
where the truncation error $R_q^n=^c\!D_q^\alpha x(t_n)-\triangle_q^\alpha x(t_n)$ is given by \rf{3.7}. Let the error function $e_n=x^n-x(t_n)$. From Eqs. \rf{5.2} and \rf{5.8}, we obtain
\begin{eqnarray*}
b_ne_n=\sum_{k=1}^{n-1}(b_{k+1}-b_k)e_k+\Gamma_q(1-\alpha)\big(f(t_n,x^n)-f(t_n,x(t_n))\big)+\Gamma_q(1-\alpha)R_q^n\,.
\end{eqnarray*}
Then, it follows a similar argument to that of Theorem 2
\begin{equation}
|e_n|\leq \Gamma_q(1-\alpha)t_n^{\alpha}\max_{1\leq k\leq n}|f(t_k,x^k)-f(t_k,x(t_k))|+\Gamma_q(1-\alpha)t_n^\alpha\max_{1\leq k\leq n}|R_q^k|,\;n\geq 1.\la{5.9}
\end{equation}
Hence, using condition \rf{5.3}, it yields
$$
|e_n|\leq L_1\max_{1\leq k\leq n}|e_k|+\Gamma_q(1-\alpha)t_n^\alpha\max_{1\leq k\leq n}|R_q^k|,\;n\geq 1.
$$
Which implies
\begin{equation}
|e_n|\leq \frac{1}{1-L_1}\Gamma_q(1-\alpha)t_n^\alpha\max_{1\leq k\leq n}|R_q^k|,\;n\geq 1.\la{5.10}
\end{equation}
Using Theorem 1 and noting that $t_k=\triangle t_k/(1-q)$, we obtain
\begin{equation}
|R_q^k|\leq \frac{1}{4\Gamma_q(1-\alpha)}\frac{1}{1-q^2}\frac{(1-q)^\alpha}{q^\alpha-q}\triangle t_k^{2-\alpha}\max_{0\leq t\leq t_k}|D_q^2x(t)|\,.\la{5.11}
\end{equation}
Substituting this into \rf{5.10}, it yields
$$
|e_n|\leq \frac{1}{1-L_1}\frac{1}{4}\frac{1}{1-q^2}\frac{(1-q)^\alpha}{q^\alpha-q}t_n^\alpha\max_{1\leq k\leq n}\triangle t_k^{2-\alpha}\max_{0\leq t\leq t_n}|D_q^2x(t)|,\;n\geq 1.
$$
The proof is completed, noting that $\displaystyle{\max_{1\leq k\leq n}}\triangle t_k=\triangle t_n,\,t_n=\triangle t_n/(1-q)$.\zb

Since $\triangle t_n=t_n-t_{n-1}=bq^{N-n}(1-q)$, then from \rf{5.7} we also obtain
\begin{equation}
|x(t_n)-x^n|\leq \frac{1}{1-L_1}\frac{1}{4}\frac{1-q}{q^\alpha-q}b^{2}q^{2(N-n)}\max_{0\leq t\leq t_n}|D_q^2x(t)|,\;\,n\geq 1.\la{5.12}
\end{equation}
This implies that for any fixed mesh point $t_n\in (0,b)$ with $n\leq (1-\delta)N$ where $0<\delta<1$, the difference solution $x^n\rightarrow x(t_n)$ as $N\rightarrow \infty$ and the convergence rate is of $O(q^{2\delta N})$-order.\\
{\bf Remark 2.}\quad Since the $q$-mesh step sizes  $\triangle t_1<\triangle t_2<\cdots<\triangle t_N$ and $\triangle t_N=(1-q)b$, therefore, in general speaking, the difference solution $x^n$ defined on the time scale $T_{q,b}$ can not admits a global convergence rate that is valid for all $q$-mesh points $\{t_n\}$.
\section*{6. NUMERICAL EXPERIMENT}
In this section, we use numerical examples to verify the high accuracy and effectiveness of the $L_{1,q}$ difference formula. In the numerical experiment, all computations are carried out by using Mathematica 11.3.\\
\commentout{
{\bf Example 1.}\quad We first test the accuracy of the $L_{1,q}$ difference formula. For this end, we take the test function:
$$
x(t)=t^2+1,\;\;\;^cD_q^\alpha x(t)=\frac{q+1}{\Gamma_q(1-\alpha)}t^{2-\alpha}B_q(1-\alpha,2)=\frac{q+1}{\Gamma_q(3-\alpha)}t^{2-\alpha},\;0\leq t\leq 1.
$$
Denote the truncation error
$$
E_q^\alpha(t_n)=|\triangle_q^\alpha x(t_n)-D_q^\alpha x(t_n)|,\;t_n=q^{N-n},\,1\leq n\leq N.
$$
Table 1 and Table 2 give the numerical results for different parameters $\alpha$ and $q$. We observe that the accuracy of the $L_{1,q}$ difference formula is very high, in particular, at the last mesh point $t_N=1$, this formula still maintains a high accuracy.
\begin{table}[htbp]
\begin{center}
{Table 1.}\quad Numerical results and errors for $q=2/3, N=10$.\\[0.15cm]
\begin{tabular}[b]{cccccccccccc}
\hline\smallskip
$\alpha$&$t_n$&$\Delta_q^\alpha x(t_n)$&$^cD_q^\alpha x(t_n)$&$E_q^\alpha(t_n)$\\
\hline\smallskip
$\alpha=0.2$&$t_2$&$3.26319\mathrm{E}{-3}$&  $3.30387\mathrm{E}{-3}$&$4.06816\mathrm{E}{-5}$ \\

&$t_5$&$2.95000\mathrm{E}{-2}$&$2.95064\mathrm{E}{-2}$&$6.37672\mathrm{E}{-6}$\\

&$t_N$&$1.13432$&$1.13432$&$5.02805\mathrm{E}{-7}$\\

$\alpha=0.6$&$t_2$&$3.09565\mathrm{E}{-2}$&  $3.15469\mathrm{E}{-2}$&$3.64605\mathrm{E}{-4}$ \\

&$t_5$&$1.35765\mathrm{E}{-1}$&$1.35796\mathrm{E}{-1}$&$2.79151\mathrm{E}{-5}$\\

&$t_N$&$1.54680$&$1.54680$&$9.29393\mathrm{E}{-7}$\\

$\alpha=0.8$&$t_2$&$1.46995\mathrm{E}{-2}$&  $1.50641\mathrm{E}{-2}$&$5.90440\mathrm{E}{-4}$ \\

&$t_5$&$8.26760\mathrm{E}{-2}$&$8.27039\mathrm{E}{-2}$&$3.05119\mathrm{E}{-5}$\\

&$t_N$&$1.41307$&$1.41307$&$6.57513\mathrm{E}{-7}$\\
\hline
\end{tabular}
\end{center}
\end{table}

\begin{table}[htbp]
\begin{center}
{Table 2.}\quad Numerical results and and errors for $q=1/2, N=10$.\\[0.15cm]
\begin{tabular}[b]{cccccccccccc}
\hline\smallskip
$\alpha$&{$t_n$}&$\Delta_q^\alpha x(t_n)$&$^cD_q^\alpha x(t_n)$&$E_q^\alpha(t_n)$\\
\hline\smallskip
$\alpha=0.2$&$t_1$&$1.46467\mathrm{E}{-5}$&  $1.39207\mathrm{E}{-5}$&$7.25979\mathrm{E}{-7}$ \\
&$t_4$&$6.18445\mathrm{E}{-4}$&$6.18415\mathrm{E}{-4}$&$2.96073\mathrm{E}{-8}$\\

&$t_7$&$2.61134\mathrm{E}{-2}$&$2.61134\mathrm{E}{-2}$&$2.30151\mathrm{E}{-9}$\\

&$t_N$&$1.10262$&$1.10262$&$1.88455\mathrm{E}{-10}$\\

$\alpha=0.6$&$t_1$&$1.74949\mathrm{E}{-4}$&  $2.11267\mathrm{E}{-4}$&$3.63180\mathrm{E}{-5}$ \\

&$t_4$&$3.88258\mathrm{E}{-3}$&$3.88291\mathrm{E}{-3}$&$3.32296\mathrm{E}{-7}$\\

&$t_7$&$7.13647\mathrm{E}{-2}$&$7.13647\mathrm{E}{-2}$&$1.08623\mathrm{E}{-8}$\\

&$t_N$&$1.31162$&$1.31162$&$3.85586\mathrm{E}{-10}$\\

$\alpha=0.8$&$t_1$&$5.95878\mathrm{E}{-4}$&  $7.91374\mathrm{E}{-4}$&$1.95496\mathrm{E}{-4}$ \\

&$t_4$&$9.59540\mathrm{E}{-3}$&$9.59599\mathrm{E}{-3}$&$5.87836\mathrm{E}{-7}$\\

&$t_7$&$1.16358\mathrm{E}{-1}$&$1.16358\mathrm{E}{-1}$&$1.24069\mathrm{E}{-8}$\\

&$t_N$&$1.41093$&$1.41093$&$2.89842\mathrm{E}{-10}$\\
\hline
\end{tabular}
\end{center}
\end{table}
}
{\noindent\bf Example 1.}\quad In this example, we use the difference scheme \rf{5.2} to solve the linear $q$-fractional differential equation:
\begin{equation}
^c D_q^{\frac{1}{2}}x(t)=\frac{1+q}{\Gamma_q(3/2)}t^\frac{3}{2}+\frac{1}{\Gamma_q(3/2)}\sqrt{t},\quad 0<t\leq 1,\,t\in T_{q,1},\,\quad x(0)=1.\la{6.1}
\end{equation}
The exact solution is $x(t)=t^2+t+1$. The numerical results are given in Table 1.

\begin{table}[h]
\begin{center}
{ Table 1.}\quad Numerical results for problem \rf{6.1}, $q=1/4, N=10$.\\[0.15cm]
\begin{tabular}[b]{cccc}
\hline\smallskip
$t_n=q^{N-n}$&$x^n$&$x(t_n)$&$|x(t_n)-x^n|$\\
\hline\smallskip
$(1/4)^{9}$&$1.0000$&$1.0000$&$3.2911\mathrm{E}{-7}$ \\
$(1/4)^{8}$&$1.0000$&$1.0000$&$1.5456\mathrm{E}{-7}$\\
$(1/4)^{7}$&$1.0001$&$1.0001$ &$7.5569\mathrm{E}{-8}$\\
$(1/4)^{6}$&$1.0002$&  $1.0002$&$2.8026\mathrm{E}{-8}$\\
$(1/4)^{5}$&$1.0010$&  $1.0010$&$1.3998\mathrm{E}{-7}$ \\
$(1/4)^{4}$&$1.0039$&  $1.0039$&$2.5337\mathrm{E}{-6}$ \\
$(1/4)^{3}$&$1.0159$&  $1.0159$&$4.0685\mathrm{E}{-5}$\\
$(1/4)^{2}$&$1.0664$&  $1.0671$&$6.5104\mathrm{E}{-4}$\\
$(1/4)^{1}$&$1.3125$&  $1.3129$&$1.4035\mathrm{E}{-4}$\\
$(1/4)^0$&$3.0000$&$3.0002$&$2.1218\mathrm{E}{-4}$\\
\hline
\end{tabular}
\end{center}
\end{table}

\noindent{\bf Example 2.}\quad In this example, we use the difference scheme \rf{5.2} to solve the nonlinear $q$-fractional differential equation:
\begin{equation}
^c D_q^{\frac{2}{3}}x(t)=\frac{1+q}{\Gamma_q(7/3)} (x-1)^{\frac{2}{3}},\quad 0<t\leq 1,\;t\in T_{q,1},\quad x(0)=1.\la{6.2}
\end{equation}
By a straightforward computation, using formulas \rf{2.4} and \rf{3.1},  we can verify that $x(t)=t^2+1$ is a solution of problem \rf{6.2}. The numerical results are given in Table 2.

From the numerical results in Table 1 and Table 2, we see that our difference method is stable, high accurate and effective.
\begin{table}[htbp]
\begin{center}
{Table 2.}\quad Numerical results for problem \rf{6.2}, $q=2/3, N=10$.
\begin{tabular}[b]{ccccc}
\hline\smallskip
$t_n=q^{N-n}$&$x(t_n)$&$x^n$&$|x(t_n)-x^n|$\\
\hline\smallskip
$(2/3)^{9}$&$1.0007$&$1.0005$&$2.0209\mathrm{E}{-4}$ \\
$(2/3)^{8}$&$1.0015$& $1.0014$&$1.6013\mathrm{E}{-4}$ \\
$(2/3)^{7}$&$1.0034$&$1.0033$&$1.3277\mathrm{E}{-4}$\\
$(2/3)^{6}$&$1.0077$&$1.0076$&$1.1208\mathrm{E}{-4}$\\
$(2/3)^{5}$&$1.0173$&$1.0172$&$9.4942\mathrm{E}{-5}$\\
$(2/3)^{4}$&$1.0390$&$1.0389$&$7.9273\mathrm{E}{-5}$\\
$(2/3)^{3}$&$1.0878$&$1.0877$&$6.2639\mathrm{E}{-5}$\\
$(2/3)^{2}$&$1.1975$&$1.1975$&$4.0754\mathrm{E}{-5}$\\
$(2/3)^1$&$1.4444$&$1.4444$&$4.6923\mathrm{E}{-6}$\\
$(2/3)^0$&$2.0000$&$2.0001$&$6.5151\mathrm{E}{-5}$\\
\hline
\end{tabular}
\end{center}
\end{table}

\section*{7. CONCLUSION }
\setcounter{section}{7} \setcounter{equation}{0}
An unconditionally stable difference formula is first presented to discretize the fractional $q$-derivative $^cD^{\alpha}_q x(t)$ of Caputo type with $0<q,\,\alpha<1$ on the time scale set $T_{q,b}$. The rigours truncation error boundness is derived.  This difference formula can be used to solve the $q$-fractional differential equation: $^cD^{\alpha}_q x(t)=f(t,x(t))$ on set $T_{q,b}$. We prove the unique existence and stability of the difference solution and derive an error estimate of $\triangle t_n^2$-order. Numerical experiments show the effectiveness and high accuracy of this difference method. This difference formula provides a useful tool for solving numerically the $q$-fractional problems on the time scale set.
\section*{ACKNOWLEDGMENTS}
The authors would like to thank the anonymous referees for many
helpful suggestions which improved the presentation of this paper.

This work was supported by the State Key Laboratory of Synthetical
Automation for Process Industries Fundamental Research Funds (Grant No. 2013ZCX02).

\end{document}